\documentclass[reprint,
 amsmath,amssymb,
 aps,twocolumn
]{revtex4}

\usepackage{graphicx}
\usepackage{dcolumn}
\usepackage{bm}
\usepackage{stmaryrd}

\DeclareMathOperator{\dist}{dist}

\begin{document}

\preprint{APS/123-QED}

\title{How densely can spheres be packed with moderate effort in high dimensions?}

\author{Veit Elser}
\affiliation{Department of Physics, Cornell University\\
Ithaca, NY 14853}


\begin{abstract}
We generate non-lattice packings of spheres in up to 22 dimensions using the geometrical constraint satisfaction algorithm RRR. Our aggregated data suggest that it is easy to double the density of Ball's lower bound, and more tentatively, that the exponential decay rate of the density can be improved relative to Minkowski's longstanding 1/2. 
\end{abstract}

\maketitle

\section{Introduction}

The packing of congruent spheres in Euclidean space has important practical implications and is a seemingly unbounded source of theoretical questions. A major recent success was the discovery by Viazovska \cite{viazovska2017sphere} and coworkers \cite{cohn2017sphere} of modular functions that make the Cohn-Elkies linear programming density upper bound \cite{cohn2003new,cohn2002new} sharp for the $E_8$ and Leech lattices, proving that these lattice-based schemes for packing spheres are the densest possible in eight and 24 dimensions. By contrast, the subject of lower bounds on achievable densities is much murkier and progress seems to have stalled.

Minkowski \cite{minkowski1905diskontinuitatsbereich} was the first to find a lower bound for general dimension $n$ that was superior to the density achieved by any of the known packing schemes available for arbitrary $n$. For example, a simple scheme is to center the spheres on the $n$-dimensional checkerboard lattice $D_n$, the subset of integer lattice points with even coordinate sum. This gives the optimal density for $n=3$ \cite{hales2005proof} and is also believed to be the best possible for $n=4$ and $5$. On the other hand, the density $\Delta$, or fraction of space covered by spheres, decays as
\begin{equation*}
\Delta(D_n)\sim \frac{1}{\sqrt{4\pi n}}\left(\frac{e\pi}{n}\right)^{n/2}\;,
\end{equation*}
which is much faster than Minkowski's bound whose leading behavior is $2^{-n}$.

Like Minkowski's result, recent advances are also based on lattices and have asymptotic densities
\begin{equation*}
c\,\frac{n}{2^n}\;,
\end{equation*}
with improvements in the value of the constant $c$. The current best bound, for general $n$, is Ball's bound \cite{ball1992lower}
\begin{equation*}
\Delta > \Delta_\mathrm{B}=\frac{n-1}{2^{n-1}}\,\zeta(n)\;,
\end{equation*}
corresponding to $c=2$. Ball's result hinges on a lemma in Bang's proof of the ``plank problem" \cite{bang1951solution} and corresponds geometrically to the transformed problem of custom-fitting a thin oblate ellipsoid in the integer lattice that avoids all lattice points except the origin. Though only the existence of the ellipsoid is established, and the corresponding packing is not explicitly constructed, the value $c=2$ appears as a sharp estimate because the ellipsoid is constrained all over its surface. Vance \cite{vance2011improved} was able to further improve $c$ when $n$ is divisible by four and Venkatesh \cite{venkatesh2013note} found that $c$ could be replaced by $\log\log n$ for very special $n$.

These increasingly sophisticated bounds, all based on lattices, stand in stark contrast to a bound that makes no reference to lattices at all and can be proved in five sentences. A set of sphere centers $S_n^*$ is ``saturated" for spheres of radius $r$ if it is impossible to add another sphere, also of radius $r$, without intersecting an existing sphere. This property implies all points not covered by a sphere are within distance $2r$ of one of the sphere centers. By doubling all the sphere radii, all of these points will be covered as well. But this could not happen if $\Delta(S_n^*)<2^{-n}$, since doubling the radii increases each sphere volume by $2^n$. We therefore know that
\begin{equation*}
\Delta(S_n^*)\ge 2^{-n}\;.
\end{equation*}
Like the lattice-based bounds, this construction is not constructive in a practical sense. On the other hand, it reveals that matching the leading asymptotic part of the sophisticated bounds is already achieved by a greedy algorithm. Information on where spheres can be placed is provided by the Voronoi diagram of the sphere centers already placed, something that can be locally updated in a sequential construction of a periodic saturated packing.

The crudeness of saturated packings suggests that easy improvements on the lower bound should be possible just by dropping the lattice constraint. Theoretically this proposal is still difficult because no one knows how to even mildly enhance the density  in a way that is also amenable to computations. Torquato and Stillinger (TS) \cite{torquato2006new,scardicchio2008estimates} have conjectured the existence of packings that have a particular limiting form of the sphere-center autocorrelation (pair distribution function) $g_2$ in high dimensions. If such packings exist, then the dominant $2^{-n}$ behavior of the density would be improved to $b^n$, with $b\approx 0.583$. The constrained optimization of $g_2$ used by TS to obtain this $b$ can be interpreted as an infinite-dimensional linear program (LP) dual to the LP used by Cohn and Elkies \cite{cohn2003new,cohn2002new,cohn2022dual} to establish upper bounds on the density. In this setting the TS-conjectured lower bound is a rigorous lower bound on the upper bounds that can be achieved with the LP method. However, this ``lower bound on the upper bounds" may well be above realizable densities if it turns out that packings with the conjectured $g_2$ do not exist. To add perspective to the Torquato-Stillinger constant $0.583$, we note that by the Kabatiansky-Levenshtein upper bound \cite{kabatiansky1978bounds}, the decay constant of the density is less than $0.661$. 

This study was prompted by the dearth of evidence that could inform the pursuit of an improved lower bound. The blindest saturated packing construction, called random sequential addition (RSA), is an obvious source of data \cite{torquato2006random}. In RSA, sphere placement is not informed by the Voronoi diagram (e.g. filling the smallest available hole), but is sampled uniformly on the available set. Accurate estimates of the saturation densities have been obtained in up to eight dimensions \cite{zhang2013precise}. These data are consistent with $b=1/2$, and that may not be surprising because this dominant behavior was proved for the related ghost-RSA construction \cite{torquato2006exactly}.

Physics inspired constructions, also called simulations, have generated data in up to 12 dimensions \cite{charbonneau2011glass}. The most widely used is the Lubachevsky–Stillinger algorithm \cite{lubachevsky1991disks}, where spheres execute Newtonian dynamics with a simple billiard-type collision rule. Initially the spheres are small and easily packed into a simulation cell with periodic boundary conditions. As the spheres redistribute themselves through collisions, they are also made to slowly grow in size. Eventually the redistribution of the spheres, into more loosely packed configurations that allow continued size-growth, slows so dramatically that the simulation cannot be continued any further. The highest densities achieved by this method, in 12 dimensions, are about 88\% higher than Ball's bound \cite{charbonneau2011glass}.

We present two kinds of data of relevance to the lower bound problem. Both are generated by the general purpose relaxed-reflect-reflect (RRR) algorithm \cite{elser2021learning}. RRR is a generalization, to general constraint sets, of the constraint satisfaction algorithm used in phase retrieval. Provided the problem at hand can be formulated as the search for a point in the intersection of sets $A$ and $B$ in some Euclidean space, and projections to these sets can be computed efficiently, the application of RRR is straightforward.

Though RRR can find some of the densest-known packings, such as Best's packing in ten dimensions \cite{best1980binary,conway1994quaternary}, our focus has been on what density improvements are possible with only moderate effort. In our first application of RRR we find packings that exactly double the density of Ball's bound while taking careful account of the effort involved. Our results extend to 22 dimensions and look like they can be continued indefinitely. In the second, more ambitious application, we do not set a density goal and instead specify a moderate level of computational effort. These results extend to 19 dimensions and the resulting densities support the existence of a bound with $b>1/2$, though less convincingly than the conclusion of the first study.

A key part of both studies is to demonstrate that the experiments have crossed into the asymptotic regime, so that the data are relevant for the lower bound question. Our handle on this is to pack spheres into the smallest possible $n$-torus, whose period equals the sphere diameter. Instead of the usual ``thermodynamic limit" of physics that tries to eliminate the effects of boundaries, we look for asymptotic behavior only with respect to the dimension $n$. Even when confined to the smallest torus, the fraction of a sphere's surface available for contacts with other spheres approaches 100\% in high dimensions. We use this property to help assess whether the asymptotic regime (in dimension) has been accessed.

\section{RRR algorithm}

The RRR algorithm \cite{elser2021learning} searches for a point $x\in A\cap B$, where $A$ and $B$ are sets in some Euclidean space. The elementary operations of the algorithm are the projectors $P_A$ and $P_B$. The point $P_A(x)$ is the element of $A$ nearest to $x$, and similarly for $B$. When a set is not convex there can be multiple nearest points, but only for $x$ in a set of measure zero. Since our computations have finite precision, the implementations of the projections always output unique points.

After selecting an initial point $x$, RRR iterates the map
\begin{equation}\label{eq:RRR}
x\mapsto x'=x+\beta\left(P_B(R_A(x))-P_A(x)\right)\;.
\end{equation}
Here $R_A(x)=2P_A(x)-x$ reflects in constraint set $A$ and $\beta$ is a parameter analogous to a time step. When $A$ and $B$ are locally affine, as they are in the sphere packing problem, RRR locally converges to fixed-point sets $X^*$ associated with solutions $x^*\in A\cap B$. Though fixed-point convergence is assured by the convexity of the affine local approximations of $A$ and $B$, RRR performs something analogous to an actual search when $A$ and $B$ are not convex. One of the sets in the sphere packing problem is nonconvex. RRR derives its name from the form it takes when \eqref{eq:RRR} is entirely expressed in terms of reflectors. For additional background on this approach to solving problems see \cite{elser2007searching,lindstrom2021survey}.

\section{$A$ and $B$ for sphere packing}

We used divide-and-concur \cite{gravel2008divide} to define $A$ and $B$, where variables are given multiple copies to make constraint projections easy. In the sphere packing problem we use ${\scriptstyle \binom{N}{2}}$ pairs of $n$-tuple (sphere center) variables. For example, $x_{ij}$ is the copy of sphere center $i$ ``that cares about" sphere $j$, and vice versa for $x_{ji}$. For all pairs $(i,j)$, the projector $P_A$ moves the centers $x_{ij}$ and $x_{ji}$ by the minimum distance to make their spheres tangent if they intersect and does nothing if these sphere-copies are already disjoint. The other projector, $P_B$, implements ``concur" by restoring equality to the copies in a distance minimizing way. The total number of (replicated) variables, $N(N-1)n$, is the dimension of the space in which RRR executes the search.

Finding distance-minimizing displacements that make two spheres tangent when they intersect (set $A$), and the concurring point nearest to $N-1$ copies (set $B$), are both easy computations. Details, including complications that arise because of the torus geometry, are provided in the appendix. Our implementation also allows the metric that defines the projections to adiabatically adjust to circumstances. Though all the spheres are identical, their equivalence under constraint projections is broken already in the initial, randomly generated configuration. This has many sphere intersections, some spheres worse off than others. Details for the general metric-update heuristic we used that addresses such inequities can also be found in the appendix.

\section{Packing the $n$-torus}

Our RRR sphere packing implementation packs $N$ unit-diameter spheres in the $n$-torus of width $w$:
\begin{equation*}
\mathcal{T}_n(w)=\mathbb{R}^n/(w\mathbb{Z})^n\;.
\end{equation*}
When extended periodically this realizes a packing of Euclidean space with density
\begin{equation*}
\Delta=N\frac{v_n(1/2)}{w^n}\;,
\end{equation*}
where $v_n(r)$ is the volume of the $n$-ball of radius $r$. To keep $N$ small in high dimensions we used the smallest possible width, or the unit torus of width $w=1$. Though our packings are in many respects random, in the unit torus each sphere will at least be tangent to $2n$ others (when the packing is extended periodically). Packing more than one sphere in the unit torus first becomes possible in four dimensions, and the optimal packing, of two spheres, is unique.

Unit diameter spheres with centers $x_1$ and $x_2$ that are tangent in Euclidean $n$-space can be packed in $\mathcal{T}_n(1)$ only if $x_1$ and $x_2$ have no coordinate differences greater than $1/2$ in absolute value. When this is not the case, then one sphere has a coset representative with a smaller coordinate-difference magnitude, giving it a smaller distance to the other sphere, in violation of the packing constraint.

\begin{figure}[t]
    \centering
    \includegraphics[width=0.45\textwidth]{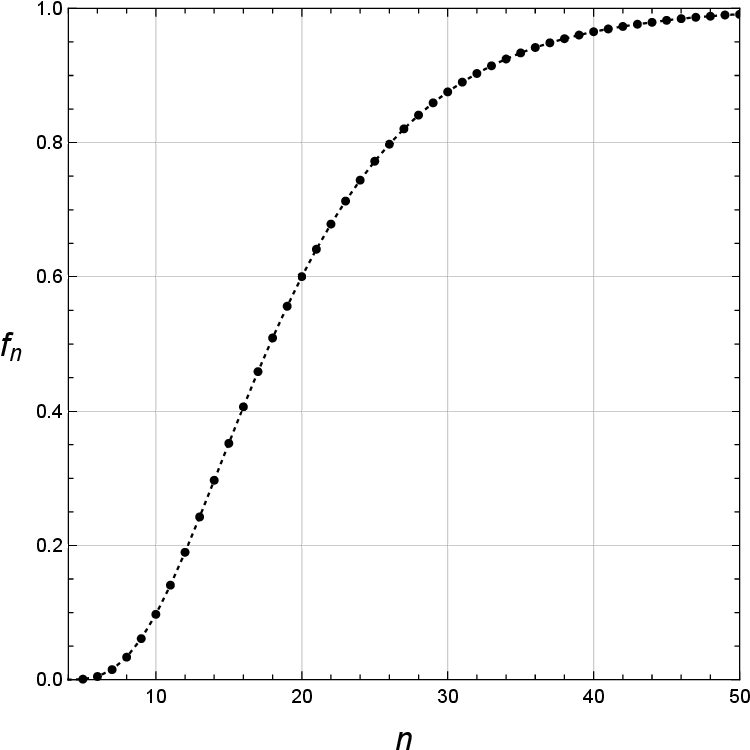}
    \caption{Fraction of the sphere of radius 1 available for tangencies when spheres are packed in the unit $n$-torus.}
    \label{fig:areafrac}
\end{figure}

The set of torus-restricted-tangencies just described, that is, the allowed set of center differences $z=x_1-x_2$, is the set
\begin{equation}\label{rsphere}
\widetilde{S}_n(1)=S_n(1)\cap C_n(1/2)\;,
\end{equation}
where $S_n(1)$ is the unit radius $n$-sphere and
\begin{equation*}
C_n(1/2)=\{z\in\mathbb{R}^n\colon \|z\|_\infty\le 1/2\}
\end{equation*}
is the centered $n$-cube of unit width.
In the appendix we provide details for the formula ($n>3$)
\begin{align}\label{rspherevol}
\mathrm{vol}(\widetilde{S}_n)&=\tilde{a}_n\nonumber\\
&=\sum_{k=0}^3 (-1)^k a_{n-k}\, b(n,k)\, c(n,k)\;,
\end{align}
where
\begin{equation*}
a_m=\mathrm{vol}(S_m(1))=m v_m(1)
\end{equation*}
is the volume of the unrestricted $m$-sphere,
\begin{equation*}
b(n,k)=\sum_{\ell=k}^3 (-1)^\ell\binom{\ell}{k} \binom{n}{\ell}
\end{equation*}
is a combinatorial factor, and
\begin{equation*}
c(n,k)=\int_{z\in C_k(1/2)}\left(1-\|z\|_2^2\right)^{\frac{n-2-k}{2}}\;,
\end{equation*}
is an integral over the $k$-cube with $c(n,0)=1$. When interpreting our packing experiments we will refer to the plot of the fraction of the tangency volumes, $f_n=\tilde{a}_n/a_n$, shown in Figure \ref{fig:areafrac}.

\section{Experiments with \textsf{toruspack}}

The main inputs for our RRR implementation, called \textsf{toruspack} \cite{toruspack}, are the dimension $n$, the number of spheres to be packed $N$, and the torus width $w$. In this study we always set $w=1$ except when we want to specify a precise value for the density. In that case we let the density determine a fractional $N$ with $w=1$, round this upward to the nearest integer, and then compensate by appropriately increasing $w$. Because $N$ grows rapidly with $n$, the value of $w$ even in these experiments is only very slightly greater than $1$.

The other inputs relate to the behavior of RRR and are not special to the sphere packing problem. Though local fixed-point convergence is fastest when $\beta=1$ is used for the RRR time step, smaller values are more productive when the search is faced with nonconvex constraints. We used $\beta=0.5$ in all the experiments. The parameter that controls the update-rate of the metric was set at the same small value, $\gamma=10^{-3}$, in all the experiments except the one instance of a hard search (Best's packing), where we used $\gamma=10^{-2}$.

\begin{figure}[t!]
    \centering
    \includegraphics[width=.45\textwidth]{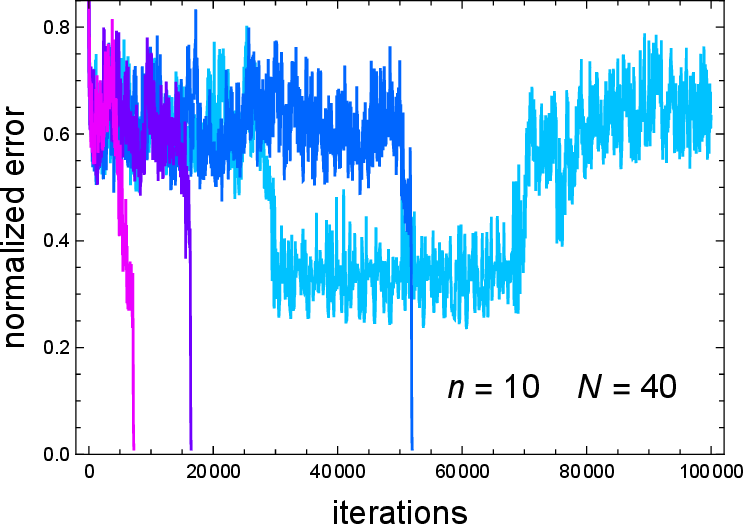}
    \caption{Error time series in four runs of packing 40 spheres in ten dimensions, three of which found Best's packing in under $10^5$ iterations.
    }
    \label{fig:besterr}
\end{figure}

The algorithm's progress is monitored via the normalized error
\begin{equation*}
\epsilon=\|x'-x\|/\sqrt{N}\;,
\end{equation*}
which may be interpreted both as the current root-mean-square 2-norm of the distances moved by the copies of each sphere, as well as the current proximity to a solution --- since $\epsilon=0$ corresponds to a fixed point. Figure \ref{fig:besterr} shows the time series of $\epsilon$ in four runs with $n=10$, $N=40$ and $w=1$. This is a hard instance, with $\epsilon$ behaving chaotically. The occasional dips to small values indicate the discovery of near solutions, where only a few concurring spheres are intersecting, or only a few spheres are unable to settle on concurred positions, or a combination of these. Eventually, when a true solution is found, $\epsilon$ goes all the way to zero. In this instance it is Best's packing \cite{best1980binary,conway1994quaternary}, the best known for $n=10$, in which the sphere centers form a binary code of Hamming distance four.

In all the experiments we terminate runs and declare them successful when $\epsilon$ falls below $10^{-4}$. The number of iterations $\mathcal{I}$ needed to find solutions is interesting because it measures the computational work, but depends unpredictably on the random initial point. Though the run-to-run variation in $\mathcal{I}$ is much smaller for our ``lower-bound instances" than it is for Best's packing, the determination of accurate values for the average number of iterations, $\overline{\mathcal{I}}$, is important for some of our results. To facilitate this, another \textsf{toruspack} input allows the user to specify the number of ``trials," or separate runs differing only in the initial random point $x_0$. The point $x_0$ is always in set $B$, where each sphere has perfectly concurring copies, but is otherwise uniformly sampled in the $n$-torus.

\begin{figure}[t!]
    \centering
    \includegraphics[width=.45\textwidth]{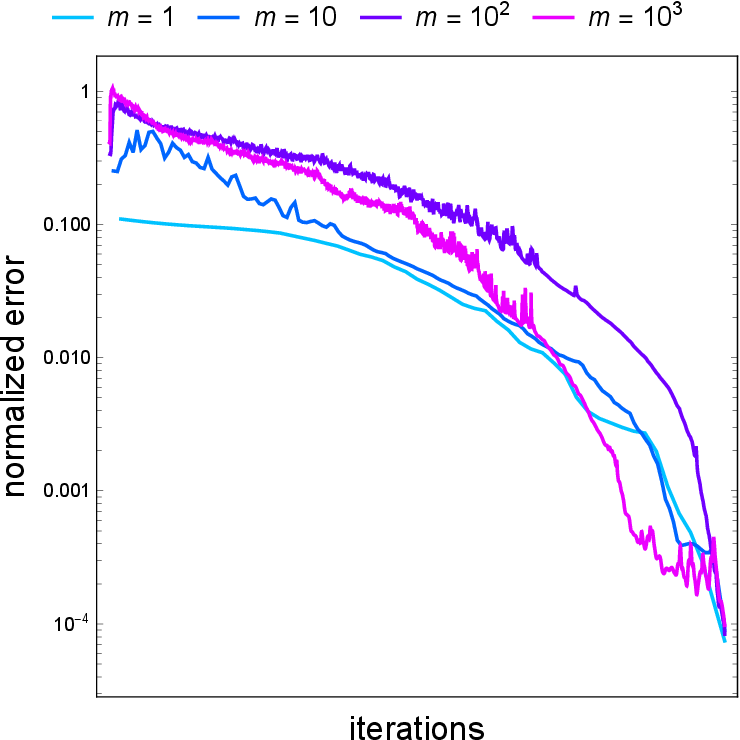}
    \caption{Representative time series of the normalized error when packing spheres in 14 dimensions as the monotonicity $m$ is increased from $1$ to $10^3$. The horizontal scales were stretched to aid comparison. Table \ref{tab:mono14} gives the actual number of iterations and the number of spheres packed for each $m$.
    }
    \label{fig:monoerr}
\end{figure}

The small run-to-run variation in $\mathcal{I}$ for the lower bound instances is related to the near monotonicity of the error time series in these relatively easy packing problems. In order to have control over the degree of difficulty, \textsf{toruspack} has a parameter $m$ called the ``monotonicity." A solution is $m$-monotone if for every iteration $i$, the error satisfies $\epsilon_{i+m}<\epsilon_i$. The easiest instances, where the error is strictly decreasing, have 1-monotone solutions. Runs are terminated by \textsf{toruspack} as soon as the $m$-monotonicity criterion is violated and the trial is declared unsuccessful. A packing instance, or $(n,N)$ pair, has $m$-monotone difficulty if a successful run with monotonicity parameter $m$ has probability one-half (when sampling initial points). In our second experiment we find the densities of packings that have 100-monotone difficulty. Using \textsf{toruspack} this means fixing $n$ and $m$, then performing trials with increasing $N$ to identify the largest $N(n,m)$ before the success probability drops below one-half. Figure \ref{fig:monoerr} shows successful error time series for packings in 14 dimensions and difficulty ranging from 1-monotone to 1000-monotone. The time scales have been stretched to have the same extent to highlight the qualitative differences. The actual number of iterations and related information is given in Table \ref{tab:mono14}. In all except the $m=1$ runs there is a very short initial rise in the error. We interpret this as transient behavior associated with the initial point exactly satisfying constraint $B$. The $m$-monotonicity criterion is imposed only after this maximum error has been passed.

The fixed point convergence shown in Figures \ref{fig:besterr} and \ref{fig:monoerr} are qualitatively different and are misleading, as presented. Not only is the overall error-behavior of the easy instances quasi-monotone, the rate of convergence accelerates. We believe this is simply related to the fact that these packings are inherently disordered, where all spheres are able to rattle around, though typically in a very small free volume. The concur projection $P_B$ at least provides a mechanism that can keep spheres from being tangent in a solution. But whereas the constraint problem for the easy instances becomes increasingly trivial with time, this is not the case in hard instances, where spheres have many tangencies in a solution. The late stage convergence of Best's packing (not shown in Figure \ref{fig:besterr}) is in fact much slower than in a disordered packing, and runs were terminated already at $\epsilon < 10^{-2}$. Best's packing, of course, was easy to confirm by rounding coordinates.

\begin{table}
\begin{center}
\renewcommand{\arraystretch}{1.2}
\setlength{\tabcolsep}{10pt}
\begin{tabular}{ l|r|r|r|r } 
 $m$ & $1$ & $10$ & $10^2$ & $10^3$\\ 
 \hline
 $N(14,m)$ & $17$ & $42$ & $115$ & $164$\\
 $\Delta/\Delta_\mathrm{B}$ & $0.392$ & $0.968$ & $2.651$ & $3.780$\\ 
 $\overline{\mathcal{I}}$ & $59$ & $163$ & $1190$ & $3980$ \\
 $\delta \mathcal{I}/\overline{\mathcal{I}}$ & $0.092$ & $0.047$ & $0.070$ & $0.184$
 \end{tabular}
\caption{Growth in the number of spheres $N(14,m)$ that can be packed in 14 dimensions as the monotonicity $m$ is increased. Also tabulated is the density relative to Ball's bound, the average number of iterations, and the normalized standard deviation in the number of iterations.
}
\label{tab:mono14}
\end{center}
\end{table}

\section{Results}

\subsection{Doubling Ball's density}

\begin{figure}[t!]
    \centering
    \includegraphics[width=.45\textwidth]{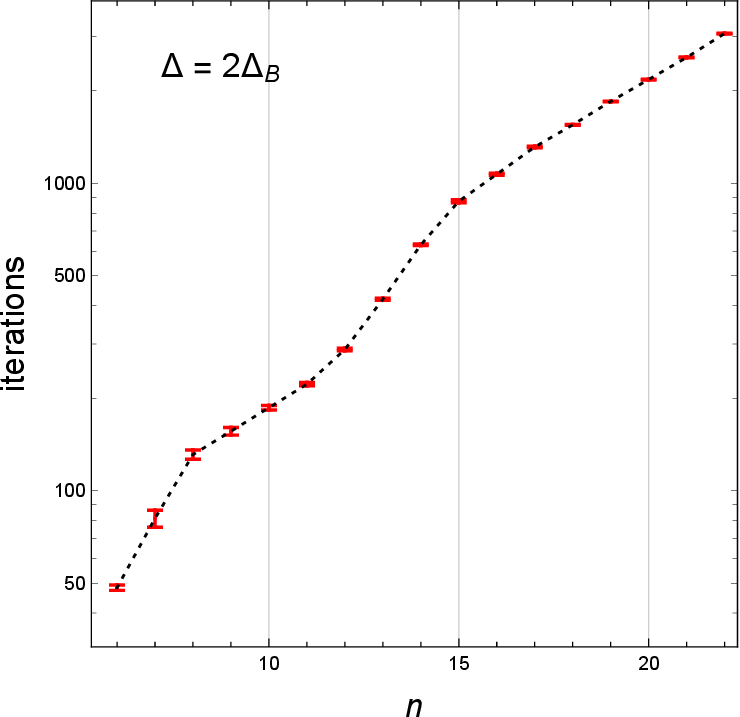}
    \caption{Average number of iterations $\overline{\mathcal{I}}_n$ used by RRR to pack spheres in $n$ dimensions at twice the density of Ball's bound.
    }
    \label{fig:doubleballiter}
\end{figure}

\begin{figure}[t!]
    \centering
    \includegraphics[width=.45\textwidth]{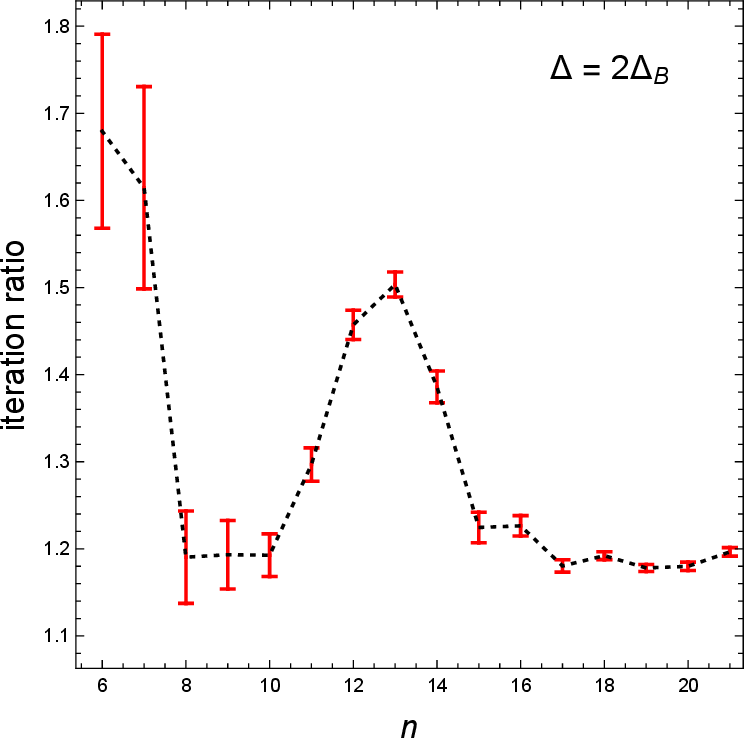}
    \caption{Successive ratios $\overline{\mathcal{I}}_{n+1}/\overline{\mathcal{I}}_{n}$ of the iteration counts plotted in Figure \ref{fig:doubleballiter}.
    }
    \label{fig:doubleballiterratio}
\end{figure}

Figure \ref{fig:doubleballiter} shows the growth in the average number of RRR iterations, $\overline{\mathcal{I}}_n$, to pack spheres at twice the density of Ball's bound. The results are consistent with simple exponential growth beyond 16 dimensions, as shown in the plot of the ratios $\overline{\mathcal{I}}_{n+1}/\overline{\mathcal{I}}_n$ in Figure \ref{fig:doubleballiterratio}. Our reach into high dimensions was limited not so much by the exponential growth in the number of iterations, but the super-exponential growth in the number of spheres being packed. For $n=22$ RRR was packing $N=11397$ spheres, and the divide-and-concur scheme works with $O(N^2)$ variables. By using neighbor lists this number can be reduced, but not by all that much since most pairs of spheres in the unit torus are neighbors for moderate $n$. However, the particulars of the growth in time and memory are mostly irrelevant for what we are aiming to demonstrate, which is showing that we can count on the algorithm being able to complete its task for arbitrary $n$. A negative result would be signs that the number of iterations might diverge at some $n$. The absence of such signs in our experiments raises confidence that Ball's lower bound on the packing density can at least be doubled.

The transient behavior in $\overline{\mathcal{I}}_n$ for $n<16$ is likely a torus artifact. The fraction $f_n$ of the unit-distance sphere available for tangencies in the $n$-torus, plotted in Figure \ref{fig:areafrac}, has an inflection point at this dimension. It appears that the torus restriction makes the task of packing spheres easier, at least for density $2\Delta_\mathrm{B}$.

\subsection{Density of 100-monotone packings}

One approach for gaining information about the possibility of $b>1/2$ in the decay of the lower bound would be to try various $b$'s and, as in the previous section, look for signs that the number of RRR iterations will diverge beyond some value of $b$. But this is unlikely to produce convincing results because one would also have to specify the equally unknown subdominant behavior of the bound, whose effect may be as large as an increase in $b$, which is surely small if nonzero.

We have taken an alternative approach, where the degree of `moderate effort' is specified through the monotonicity $m$ of the solution process. As shown in Figure \ref{fig:monoerr}, the time series of the normalized error is suggestive of a reliable route to packings for $m$ as large as $10^3$. We chose $m=100$ because this lowers the average number of iterations and also the number of spheres to be packed. The small run-to-run variation in the number of iterations, $\delta \mathcal{I}/\overline{\mathcal{I}}=0.07$ for $n=14$ (Table \ref{tab:mono14}), also lends support to the assertion that RRR can be counted on to complete its task at this degree of effort.

\begin{table}
\renewcommand{\arraystretch}{1.2}
\setlength{\tabcolsep}{10pt}
\begin{tabular}{ r|r|r } 
 $n$ & $N(n,100)$ & $\overline{\mathcal{I}}$ \\ 
 \hline
 $6$ & $4$ & $71$  \\
 $7$ & $5$ & $237$ \\ 
 $8$ & $8$ & $293$  \\
 $9$ & $12$ & $454$  \\
 $10$ & $18$ & $471$  \\ 
 $11$ & $29$ & $588$  \\ 
 $12$ & $46$ & $736$  \\ 
 $13$ & $72$ & $869$  \\ 
 $14$ & $115\pm 2$ & $1190$  \\ 
 $15$ & $156\pm 1$ & $1080$  \\ 
 $16$ & $274\pm 2$ & $1410$  \\ 
 $17$ & $507\pm 1$ & $1960$  \\ 
 $18$ & $967\pm 2$ & $2830$  \\ 
 $19$ & $1884\pm 3$ & $4240$  
\end{tabular}
\caption{Number of spheres $N(n,100)$ that RRR can pack in the unit $n$-torus with 100-monotone difficulty, along with the average number of iterations. These are the data that were used for the density shown in Figure \ref{fig:monodens}.
    }
    \label{tab:monodens}
\end{table}

The number of spheres $N(n,100)$ that can be packed by the $m=100$ criterion are given in Table \ref{tab:monodens}. For small $n$ these numbers are known with high accuracy. That is because enough trials can be performed to establish that $N(n,100)$ spheres have 100-monotone solutions in over 50\% of trials, while this drops below 50\% for $N(n,100)+1$ spheres. For larger $n$ the uncertainties in $N(n,100)$ are estimated as follows. $T$ trials are performed for a range of equally spaced candidate values of $N(n,100)$ and the number of successful 100-monotone solutions $S$ is tabulated for each of them. The success probability is estimated as $p=S/T$, with standard uncertainty $\sigma=\sqrt{p(1-p)/T}$. The smallest candidate is identified as the largest $N(n,100)$ for which $p-\sigma>1/2$, while the largest candidate is the largest $N(n,100)$ for which $p+\sigma>1/2$. The uncertainties given in Table \ref{tab:monodens} correspond to half the difference of these extreme estimates of $N(n,100)$. We performed 100 trials for each $N(n,100)$.

\begin{figure}
    \centering
    \includegraphics[width=.45\textwidth]{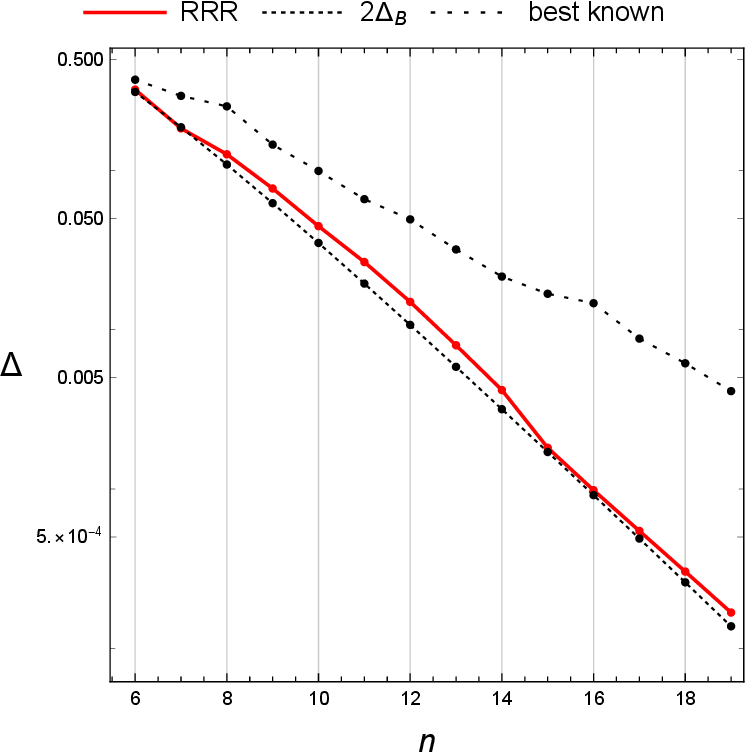}
    \caption{Density of 100-monotone packings found by RRR compared with twice Ball's bound and the best known packings.
    }
    \label{fig:monodens}
\end{figure}

\begin{figure}
    \centering
    \includegraphics[width=.45\textwidth]{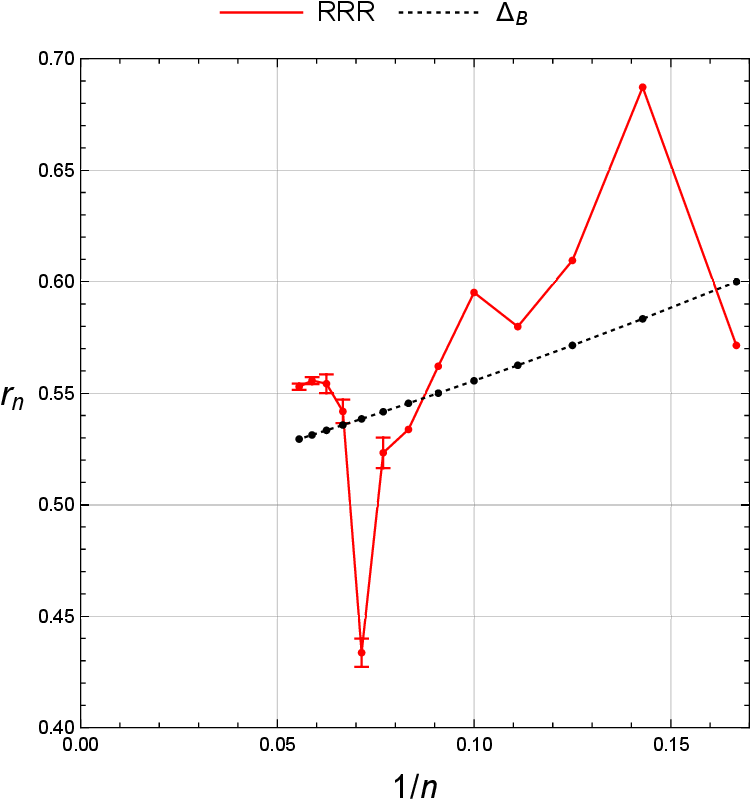}
    \caption{Successive ratios $r_n=\Delta_{n+1}/\Delta_n$ of the 100-monotone densities of Figure \ref{fig:monodens} plotted against $1/n$.
    }
    \label{fig:monodensratio}
\end{figure}

The densities associated with the data in Table \ref{tab:monodens} are plotted in Figure \ref{fig:monodens} and compared with $2\Delta_\mathrm{B}$ and the densest known packings. In this plot the uncertainties are smaller than the plotting symbol. The former become noticeable in the plot of the density ratios in Figure \ref{fig:monodensratio}. Both plots are consistent with the interpretation given earlier, that the tight unit-torus constraint makes packing easier, but that this effect wears off above 16 dimensions. Like the $2\Delta_\mathrm{B}$ experiments, here it appears we have also just barely been able to access the asymptotic regime. Though the ratios are plotted against $1/n$, the data do not extend to large enough $n$ and the uncertainties are too large to attempt an estimate of $b$ by extrapolating to $n=\infty$. Still, it appears that a lower bound with $b>1/2$ is not easily ruled out by these results.

\begin{figure}[t!]
    \centering
    \includegraphics[width=.45\textwidth]{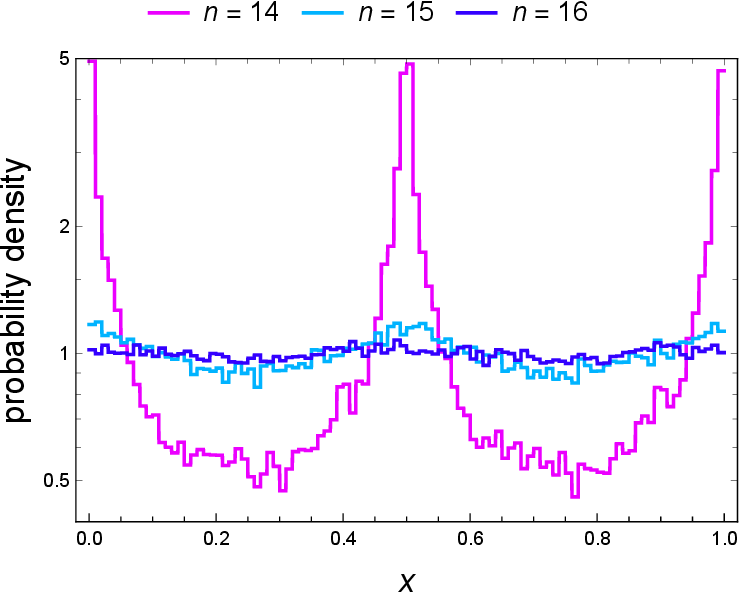}
    \caption{Disappearance of the binary-code coordinate modulation above 14 dimensions.
    }
    \label{fig:coordDistribution}
\end{figure}

\subsection{Binary-code coordinate modulation}

Additional evidence that torus artifacts are absent above $n=16$ can be seen in the distribution of coordinate values of the packed spheres, shown in Figure \ref{fig:coordDistribution}. In a disordered packing, without periodicity imposed by an $n$-torus, the coordinates should have a uniform distribution. While this is what we see above 16 dimensions, there are strong departures in lower dimensions. As epitomized by Best's packing, the unit-torus is the perfect scaffolding for packings with a Hamming-distance-4 binary code structure. Though the 100-monotone packings are disordered and far from a binary code, the distribution of coordinates shows a strong binary modulation in low dimensions. Figure \ref{fig:coordDistribution} shows how the amplitude of this modulation rapidly decays above 14 dimensions. To create these distributions we applied independent shifts $s$ to the coordinates in each of the $n$ dimensions that maximize the sums
\begin{equation*}
\sum_{i=1}^N \cos{4\pi(x_i+s)}\;.
\end{equation*}
After alignment, the distributions for all $n$ dimensions, and all the successful 100-monotone packings, were combined into a single distribution.

\subsection{Sphere center autocorrelation}

\begin{figure}[t!]
    \centering
    \includegraphics[width=.45\textwidth]{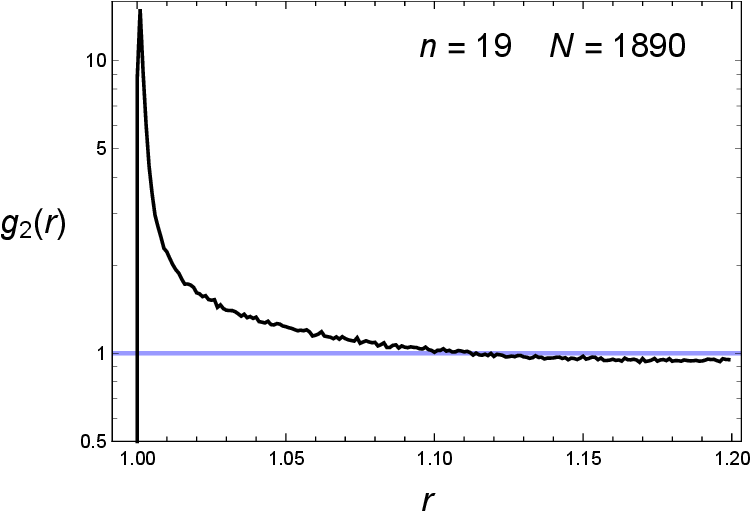}
    \caption{Sphere-center autocorrelation function $g_2(r)$ for the 100-monotone packings in 19 dimensions. 
    }
    \label{fig:g2}
\end{figure}

Because of the compact nature of our torus constraint, the sphere center autocorrelations $g_2$ of our packings are probably not good tests of the Torquato-Stillinger conjectures \cite{torquato2006new}. The $g_2(r)$ function for our highest-dimension ($n=19$) and densest ($N=1890$) packing is plotted in Figure \ref{fig:g2}. To avoid the normalization complications arising from the torus, these were created by filling center-center distance bins ($r$), first for the spheres of the packing, and also for an equal number of uniformly distributed random points. The plot shows the ratio of the two kinds of bin counts.

\section{Outlook}

The endurance of Minkowski's $2^{-n}$ leading-order behavior of the sphere packing density lower bound can be interpreted in two ways. Either it represents something fundamental that has resisted proof, mostly for lack of imagination on how that can be done. The other possibility is that packings with a slower decay exist, but that constructions analogous but better than saturated packings have likewise eluded the imagination. Though experiment will never replace proof (and imagination), it can provide guidance on which of these alternatives is the more promising one to pursue.

In this study we showed that the RRR algorithm is a good source of data for the lower bound question. Like random sequential addition and the Lubachevsky–Stillinger algorithm \cite{lubachevsky1991disks}, RRR is a simple dynamical system, appropriate for constructing the disordered packings that prevail at the lower bound. Unlike random sequential addition, but like Lubachevsky–Stillinger, RRR acts synchronously on the spheres in the packing. However, unlike Lubachevsky–Stillinger, RRR is not prone to jamming and glassy behavior.

Our two sets of results fall short of the level one has come to expect, say in statistical mechanics, for supporting hypotheses in intractable problems. Though the experiments have nearly doubled the reach into high dimensions, even higher dimensions are needed to convincingly eliminate ``torus artifacts."

Using neighbor lists and allocating divide-and-concur copies only for neighbors will certainly cut down on memory as well as the time per iteration. But only a few extra dimensions can be gained in this way. A more significant improvement is possible in principle by modifying the periodicity of the packings. Instead of using the tightest possible ``cubic" cell for unit diameter spheres, $\mathbb{Z}^n$, we propose using the checkerboard lattice $D_n$ to define the periodicity. This is not the tightest possible cell with that symmetry --- that would be $D_n/\sqrt{2}$ --- but the one having the same volume, up to a factor of two, as $\mathbb{Z}^n$. The main advantage of $D_n$ over $\mathbb{Z}^n$ is that the spheres will have a much rounder environment. For example, the analogous inflection point in the fractional volume available for tangencies (Figure \ref{fig:areafrac}) occurs already at $n=12$. This lowering of the asymptotic regime will have a greater benefit than extending the range in $n$ through neighbor lists. Also, with this scale for $D_n$ the cosets of $(\mathbb{Z}/2)^n/D_n$ support dense packings in much the same way as $(\mathbb{Z}/2)^n/\mathbb{Z}^n$ (Hamming-distance-4 codes). The disappearance of a coordinate modulation with $n$ (Figure \ref{fig:coordDistribution}) could again be used to assess the homogeneity of the packings. Hard sphere liquid simulations have also used $D_n$ boundary conditions \cite{charbonneau2022dimensional}.


\section*{Acknowledgements}

I thank Keith Ball for providing unpublished material on his lower bound.


\bibliography{refs}

\appendix

\section{Software}

The RRR C-language implementation \textsf{toruspack}, a short user's guide, and sample outputs can be found at \cite{toruspack}. The program is a single file and compiles with just the standard libraries.

\section{RRR on the torus}

Though RRR is strictly defined for constraint sets in Euclidean space, only two small modifications are required to preserve the fixed-point behavior on the torus.

In our divide-and-concur formulation of the sphere packing problem, the variables live in the product space of $N(N-1)$ tori $\mathcal{T}_n(w)$, and it is enough to address the modifications in just one such factor.
After replicating the sets $A$ and $B$ in all of $\mathbb{R}^n$ with the periodicity of $\mathcal{T}_n(w)$, the corresponding torus-distance between points $x,x'\in \mathbb{R}^n$ is defined by
\begin{equation}\label{eq:torusdist}
\dist(x,x')=\|\;[ x-x']_w\|_2\;,
\end{equation}
where $\|\;\|_2$ is the standard norm in Euclidean space and $[y]_w=y'$ is the translate of $y$ by an element of $(w\mathbb{Z})^n$ where all components of $y'$ have absolute value bounded by $w/2$. We use the notation $\|y'\|_\infty\le w/2$ to express this property. Points with some components exactly equal to $\pm w/2$, for which $[ \;]_w$ is ambiguous, do not arise in floating point computations.

The computations of $P_A$ and $P_B$ on the torus need not be modified as long as \eqref{eq:torusdist} is being minimized and the projection outputs are interpreted as coset elements. The first modification of RRR is the definition of the reflector, here for set $A$:
\begin{equation*}
R_A(x)=x+2 [ P_A(x)-x]_w\;.
\end{equation*}
In words, this corresponds to ``reflect in the nearest coset element of $A$."

The second modification is the rule for incrementing the current $x$ by the difference of projections:
\begin{equation*}
x\mapsto x'=x+\beta [ P_B(R_A(x))-P_A(x) ]_w\;.
\end{equation*}
With this modification, the local analysis of RRR in $\mathbb{R}^n$, in the presence of a set intersection or near-intersection, also applies in the torus. In particular, convergence to solutions is still characterized by the vanishing of the Euclidean error, $\|x-x'\|_2$.

\section{Disjoint-sphere projection}

Let $x_1=x_{ij}$ and $x_2=x_{ji}$ be the sphere-center copies that implement the disjointedness constraint between spheres $i$ and $j$, both of radius $r$. For packings in the torus $\mathcal{T}_n(w)$ we reduce the difference vector $u=x_1-x_2$ so that $\|u\|_\infty\le w/2$. If this reduced vector satisfies $\|u\|_2\ge 2r$, the disjointedness constraint is satisfied and the projection leaves $x_1$ and $x_2$ unchanged.

Now suppose $\|u\|_2< 2r$. To satisfy the disjointedness constraint we seek a difference vector $u'$ satisfying $\|u'\|_\infty\le w/2$, $\|u'\|_2\ge 2r$ and which minimizes $\|u'-u\|_2$. The last two conditions imply the new centers
\begin{align*}
x'_1&=x_1+(u'-u)/2\\
x'_2&=x_2-(u'-u)/2
\end{align*}
satisfy the disjointedness constraint in a distance minimizing way. If we did not have the first constraint, $\|u'\|_\infty\le w/2$, then $u'$ would just be a rescaling of $u$:
\begin{equation*}
\tilde{u}=\frac{2r}{\|u\|_2}u
\end{equation*}
To describe what we must do when at least one component of $\tilde{u}$ exceeds $w/2$ in absolute value, first express $u$ as the orthogonal decomposition
\begin{equation*}
u=u_\parallel+u_\perp\;,
\end{equation*}
where $\parallel$ corresponds to all the components of $\tilde{u}$ whose absolute values exceed $w/2$. In this decomposition the distance minimizing difference vector has the form
\begin{equation*}
u'=(w/2)\,\mathrm{sgn}(u_\parallel)+u'_\perp\;,
\end{equation*}
where the sign function $\mathrm{sgn}(\;)$ is defined to be zero on the $\perp$ components. The optimization problem for $u'$ has been reduced to a similar optimization problem, now for the vector $u'_\perp$. This vector has fewer components, must also satisfy $\|u'_\perp\|_\infty\le w/2$, and has a reduced magnitude bound:
\begin{equation}\label{eq:uperpbound}
\|u'_\perp\|_2^2\ge (2r)^2-(w/2)^2\|\mathrm{sgn}(u_\parallel)\|_2^2\;.
\end{equation}
This recursive definition of $u'$ terminates when all components of the rescaled vector $\tilde{u}_\perp$ have absolute value below $w/2$. The level of recursion never exceeds four, even in the smallest torus ($w=2r$), since in that case the bound in \eqref{eq:uperpbound} collapses to zero when $u_\parallel$ has four nonzero components.

\section{Concur projection}

In the concur projection, the sphere-center copies $x_{i1},x_{i2},\ldots,x_{iN}$ are all replaced by the same point in a distance minimizing way. Without the complication of the torus geometry, the distance minimizing point is just the centroid. To simplify the presentation, when taking account of the torus, we may treat the $n$ dimensions independently. Given coordinates $y_1,\ldots,y_N\in\mathcal{T}_1(w)=\mathbb{R}/(w\mathbb{Z})$ (in one of the dimensions), the problem is to find the $\bar{y}\in \mathcal{T}_1(w)$ that minimizes
\begin{equation}\label{eq:concurdist}
\sum_{j=1}^N\; [ \bar{y}-y_j]_w^2\;.
\end{equation}

For any $\bar{y}$, the numbers $[ \bar{y}-y_j]_w$ are ordered in the interval $B_w=[-w/2,w/2]$, which extends to an ordering on the circle if we identify $-w/2$ with $w/2$. Different $\bar{y}$ correspond to different ways of cutting the circle. There are exactly $N$ places to cut the circle, each case corresponding to a different choice of $j$ such that $[ \bar{y}-y_j]_w$ is the smallest element in $B_w$. The ordering on the circle determines the ordering of all the others, now in $B_w$, and the distance minimizing $\bar{y}$ is simply the centroid of the correspondingly ordered $y_j$. The $N$ centroids $\bar{y}$ determined by the $N$ places to cut the circle will have $N$ squared-distances \eqref{eq:concurdist}, and the projection selects among these the smallest.

Instead of computing $N$ centroids and their associated squared-distances, the projection can be computed more efficiently by calculating the changes in the squared distance between successive cuts of the circle. The work in keeping track of the smallest squared-distance after all $N$ changes have been computed scales as $O(N)$, no different from the centroid computation without the torus complication.

\section{Metric auto-tuning}

A diagonal modification of the isotropic Euclidean metric preserves the local convergence of the RRR algorithm. In the sphere packing problem, first without the torus complication, a natural diagonal modification is
\begin{equation*}
\dist^2(x,x')=\sum_{(i,j)} g_{ij}\left(\|x_{ij}-x'_{ij}\|_2^2+\|x_{ji}-x'_{ji}\|_2^2\right)\;,
\end{equation*}
where the sum is over all pairs of sphere centers and the $g_{ij}$ are positive parameters. Because a rescaling of the variables (and their constraints) restores the Euclidean metric, local convergence (of the rescaled variables) still holds. We will keep the variables (and constraints) unrescaled and adjust the metric parameters $g_{ij}$ by a heuristic based on constraint discrepancies. The adjustment is performed automatically, and adiabatically, so as not to upset the local convergence.

The local constraint discrepancy is defined as
\begin{equation*}
\epsilon^2_{ij}=\|x^A_{ij}-x^B_{ij}\|_2^2+\|x^A_{ji}-x^B_{ji}\|_2^2\;,
\end{equation*}
where superscripts $A$ and $B$ denote the $A$ and $B$ projections in the current RRR update. The sphere-center copies $x_{ij}$ and $x_{ji}$ for which $\epsilon^2_{ij}$ is above average are the most troublesome, and deserve an above-average metric weight. We implement this heuristic with the following metric-weight update rule, in each RRR iteration,
\begin{equation*}
g'_{ij}=g_{ij}+\gamma\left(\epsilon^2_{ij}/\langle\epsilon^2\rangle-g_{ij}\right)
\end{equation*}
where $\langle\;\rangle$ denotes the average of $\epsilon^2_{ij}$ over $(i,j)$ and $\gamma>0$ is a small parameter. The metric parameters are initialized at $g_{ij}=1$ and over time the most troublesome constraints have their weights increased. As the packing is refined, the metric parameters approach an equilibrium distribution $g_{ij}\approx \epsilon^2_{ij}/\langle\epsilon^2\rangle$ over a characteristic time of $1/\gamma$ iterations.

The metric parameters have no effect on the projection to the $A$ constraint, since the pairs of sphere centers involved share the same weight. However, the weights appear in an intuitive way in the $B$ constraint. The expression \eqref{eq:concurdist} being minimized for concurrence of the copies is replaced by
\begin{equation*}
\sum_{j=1}^N\, g_{i j}\, [ \bar{y}-y_j]_w^2\;,
\end{equation*}
which has the effect of weighting copy $y_j=x_{ij}$ of sphere center $i$ with weight $g_{ij}$ in the centroid computations. The copy in the most troublesome constraint thereby gets a stronger voice in deciding the sphere's position.

\section{Torus-restricted sphere volume}

The volume $\tilde{a}_n$ of the set $\widetilde{S}_n(1)$ defined in \eqref{rsphere} can be expressed as
\begin{equation}\label{vol1}
\tilde{a}_n=\prod_{i=1}^n\left(\int_{-\infty}^\infty dx_i\, \theta(x_i)\right) 2\delta\!\left(\|x\|_2^2-1\right)\;,
\end{equation}
where $\delta(\;)$ is the Dirac delta distribution and
\begin{equation*}
\theta(x)=\left\{
\begin{array}{cc}
1, & |x|\le 1/2\\
0, & \mbox{otherwise.}
\end{array}
\right.
\end{equation*}
Rewriting \eqref{vol1} in terms of $\bar{\theta}=1-\theta$, we notice that terms involving products of four or more $\bar{\theta}$ are zero because $\|x\|_2^2 > 1$ when four of the $x_i$ have absolute value greater than $1/2$. The surviving terms can be indexed by $\ell$, \mbox{$0\le\ell\le 3$}, and have the same value that only depends on $\ell$ (in addition to $n$) :
\begin{equation*}
\tilde{a}_n=\sum_{\ell=0}^3(-1)^\ell\binom{n}{\ell}\, d(n,\ell)\;.
\end{equation*}
The integral in
\begin{multline}\label{dnl}
d(n,\ell)=\prod_{j=0}^\ell\left(\int_{-\infty}^\infty dy_j\, (1-\theta(y_j))\right)\\
\times\int_0^\infty  (a_{n-\ell}\,r^{n-\ell-1} dr)\, 2\delta\!\left(\|y\|_2^2+r^2-1\right)
\end{multline}
over the $n-\ell$ variables that only appear via their norm $r$ in the integrand has been rewritten in terms of the volume $a_{n-\ell}$ of the sphere in that many dimensions. The final step of the derivation of formula \eqref{rspherevol} is to express \eqref{dnl} as a sum over the number of factors of $\theta$, indexed by $k$, all of which are equal:
\begin{equation*}
d(n,\ell)=\sum_{k=0}^3(-1)^k \binom{\ell}{k}\, e(n,k)\;.
\end{equation*}
Expressing the integral in $e(n,k)$ as earlier for $d(n,\ell)$,
\begin{multline*}
e(n,k)=\prod_{j=0}^k\left(\int_{-\infty}^\infty dy_j\, \theta(y_j)\right)\\
\times\int_0^\infty  (a_{n-k}\,r^{n-k-1} dr)\, 2\delta\!\left(\|y\|_2^2+r^2-1\right)\;,
\end{multline*}
we obtain
\begin{align*}
e(n,k)&=a_{n-k}\prod_{j=0}^k\left(\int_{-\infty}^\infty dy_j\, \theta(y_j)\right)(1-\|y\|_2^2)^{\frac{n-k-2}{2}}\\
&=a_{n-k}\, c(n,k)\;.
\end{align*}

\clearpage

\end{document}